\documentclass[11pt]{article}
\usepackage{latexsym}
\usepackage[dvips]{graphics}

\newfam\msbfam
\font\twelmsb=msbm10 at 12pt

\font\sevenmsb=msbm10 at 7pt
\font\fivemsb=msbm10 at 5pt

\textfont\msbfam=\twelmsb
\scriptfont\msbfam=\sevenmsb

\scriptscriptfont\msbfam=\fivemsb

\def\Bbb{\fam\msbfam\twelmsb}

\def\C{{\Bbb C}}
\newcommand{\R}{{\Bbb R}}
\newcommand{\ar}{\longrightarrow}
\def\p{{\Bbb P} }

\newcommand{\z}{\zeta}

\newcommand{\T}{{\Bbb T}}
\newcommand{\gs}{{\rm H}^0(\p_1, {\cal O}(\T))}

\newcommand{\rl}{{\rm L}}

\newcommand{\Th}{\theta}
\newcommand{\om}{\omega}
\newcommand{\Z}{{\Bbb Z}}
\newcommand{\Om}{\Omega}
\newcommand{\di}{\displaystyle}

\newcommand{\be}{\beta}
\newcommand{\un}{\underline}
\newcommand{\ff}{{\rm f}}

\begin{document}

\begin{center}
{\Large\bf  Minimal Surfaces from Monopoles}

{\large Anthony  Small\footnote{On leave at: Faculty of Mathematics, University
of Southampton, Southampton SO17 1BJ, England.}
 
 Mathematics Department

 National University of Ireland, Maynooth

 Co. Kildare, Ireland.

asmall@maths.may.ie}
 
 \end{center}
\begin{abstract}
{\small  \noindent 
The geometry of minimal surfaces generated by
charge 2  Bogomolny monopoles on $\R^3$ is described in terms of the
moduli space parameter $k$. We find that the distribution of Gaussian curvature
on the surface reflects the monopole structure. This is elucidated by the 
behaviour of the Gauss maps of the minimal surfaces.

\vspace{0.1in}
}  \end{abstract}

\vspace{0.2in}

{\bf \S 1. Introduction.}

In \cite{h1}, it was shown that the data comprising a static SU(2) monopole on
$\R^3$ are encoded in its {\it spectral curve}, an auxiliary algebraic curve in
$\T$, the total space of the holomorphic tangent bundle of $\p_1$. $\T$ is
viewed as the space of all oriented lines in $\R^3$, and the spectral curve
parameterises the monopole's {\it spectral lines}. These lines should be
thought of as  going through the locations of the monopole particles,
\cite{hmm}. For a monopole of charge $\ell$, the spectral curve is an
$\ell$-fold branched cover of $\p_1$, of virtual genus $(\ell-1)^2$. For
further details and background information see \cite{ah}, \cite{jt},
\cite{sut} and \cite{WW}.
 
 Recall that $\T$ compactifies, by the addition of
a single point at infinity, to a quadric cone ${\cal C}(Q)$, with vertex $v$
say, in $\p_3$. Classical osculation duality gives a correspondence
between full curves on ${\cal C}(Q)$, and  certain  full curves in $\p_3^{\,
\ast}=\C^3\cup v^\ast$. The latter  give null curves in $\C^3$, and hence
project to (branched) minimal surfaces in $\R^3$: in fact all 
non-planar minimal surfaces in $\R^3$ arise in this way. This correspondence
was discovered by Lie, see \cite{da}, \cite{h1}. The main features of a
minimal surface, i.e. the total Gaussian curvature, end structure, branch
points and symmetries may be read off the auxiliary curve in $\T$, see
\cite{sm3} for further details. 

Accordingly, a monopole on $\R^3$ generates, and is determined by, an
auxiliary minimal surface in $\R^3$: this was observed in \cite{a} and
\cite{h1}. (However, the correspondence is not understood directly, but via
the spectral curve.)  Two simple questions arise:

$\bullet$\ what do these minimal surfaces look like?

$\bullet$\ How does the geometry of the auxiliary surface reflect a
monopole's structure?

Unfortunately, for charge  $\ell\geq 2$, the surfaces are extremely
complicated.  However, we show that for $\ell=2$, the  key features are
tractible, describe them, and indicate briefly how they relate to the
monopole. 

Recall that the orbits of the Atiyah-Hitchin manifold ${\cal M}_2^0$, are parameterised 
by the elliptic modulus $k\in[0,1)$, see \cite{ah}. For $k\neq0$, each orbit
contains a reduced, centred monopole with  elliptic spectral curve $S_k$, as
described in $\S 2$. The main technical contribution here, which is
given in Theorem 4.7, is the derivation of a set of tractible formulae,
involving elliptic functions,  for the components of the auxiliary null curves
determined by osculation of $S_k, \ k\in(0,1)$. 

Now, for a generic  monopole of charge $\ell\in{\Bbb N}$,   the
auxiliary minimal surface is complete, finitely branched and of finite total
curvature $-4\pi\ell$,  $\ell$ being equal to the degree of the Gauss map of
the minimal surface. The global structure of such objects is well understood,
largely as a result of the seminal work of Osserman \cite{os},  see
\cite{lawson} and \cite{osser},  for further details. Their geometry is
`concentrated' in the sense that, `from infinity one sees a finite number of
planes passing through the origin', see \cite{Meeks1} for a precise statement.
 
Observe then that the total Gaussian curvature on the minimal surface equals 
   minus the total energy of the monopole. The results described here show
moreover that  for $\ell=2$, the {\it local}  distribution of Gaussian
curvature on the minimal surfaces, as $k$ varies, reflects the `monopole
dynamics'. This is closely related to the behaviour of the area measure
induced on the spectral curve by the monopole's `Gauss map', i.e. the branched
covering map $\pi:S_k\ar\p_1$.  Thus the energy of the monopole  is  tied to
the twisting of its spectral lines.

 For $\ell=2$, $k\neq0$, the auxiliary minimal surfaces are two ended Klein 
bottles, where the ends are perpendicular to the spectral lines through the
origin. The geometry of each of these surfaces is organised by the
configuration of  six branch points in the metric on the surface. These are
connected by a `pointed star' on the surface, $\Gamma_{S_k}$, formed from the
image of four of the quarter-period circles of $S_k$. The Gaussian curvature
localises on parts of this star to a degree that varies with $k$. In
particular, as $k\ar1$, and the monopoles become well-separated,
the Gaussian curvature  localises at the two monopole particles. Moreover,
in this limit, the normal lines to the minimal
surface in the vicinity of the particles become exponentially close, relative
to separation distance, to monopole spectral lines. Monopole scattering, cf.
\cite{ah},  is reflected by an exchange of Gaussian curvature in which
`particles of curvature' become more or less attenuated according to their
role in the interaction: the behaviour of the Gauss maps of the minimal
surfaces on $\Gamma_{S_k}$, as $k$ varies, elucidates this.

The paper is organised as follows. In $\S$2  we
briefly review background material about osculation duality.
In $\S 3$ we see what can be said about the charge 2 case prior to writing
down explicit formulae. $\S 4$ contains the main technical results of the
paper. These enable us to write down explicit formulae for the minimal
surfaces which elucidate a number of subtle features of their geometry. In
particular, they furnish a useful formula for the Gauss maps of these surfaces
in terms of the moduli parameter $k$.  This allows us to study the area
measure induced by the Gauss maps on the spectral curves in the limits
$k\ar0,\ 1$.    These and other issues are explored in $\S5-\S8$.


\vspace{0.3in}  

  {\bf \S 2. Osculation duality.}   

 Let  $\pi:\T\ar\p_1$ be   the projection map. Let $\z$ be an affine
coordinate on $\p_1$ and $(\z, \eta)$ be the coordinates given by  $\di{(\z,
\eta)  \ar \eta\frac{d}{d\z}}$.  $\gs\cong\C^3$ and the real structure 
$\tau:\T\ar\T$, given in  local coordinates by $(\z,
\eta)\ar(-\overline{\z}^{-1}, -\overline{\eta}\overline{\z}^{-2})$, 
determines the $\R^3$ of $\tau$-invariant {\it real sections}   of the form  
$$
\di{\sigma_{(x_1,x_2,x_3)}(\z)=((x_1+ix_2)-2x_3\z-(x_1-ix_2)\z^2)\frac{d}{d\z}}.
$$ 
 $\T$ parameterizes the oriented lines in $\R^3$  and, equivalently,  the
affine null planes in $\C^3$.   

 Intrinsically,  classical
osculation  duality may be understood in terms of  the family of global
sections  osculating a curve $S\subset\T$; this family   determines a null
holomorphic  curve  $\Omega:S^\ast\ar\C^3\cong \gs$, where $S^\ast$ is the
desingularization of $S$,  punctured at the finite number of points  that
correspond to points where $S$ osculates a fibre.   By duality, the original
curve parameterizes   the set of affine null hyperplanes  osculating  the null
curve in $\C^3$.

If $S\subset\T$, an $\ell$-fold branched covering of $\p_1$, $\ell\geq2$, is an
irreducible algebraic curve, then osculation
duality determines a (finitely) branched minimal immersion $\phi={\rm
Re}(\Omega):S^\ast\ar\R^3$. The branched metric, $ds^2$, induced on $S^\ast$,
is complete in the sense that every curve that approaches a puncture has
infinite length. (The punctures correspond   to the ends of the minimal
surface.)  However it will in general be (finitely) branched, i.e. have
isolated zeros. In summary, cf. \cite{sm3}, the geometry of the minimal
surface may be discerned from $S$ as follows:

$\bullet$ the
ends of the minimal surface correspond to the points where $S$ osculates a
fibre;

$\bullet$ the zeros of the branched metric are  caused by
hyperosculating sections;

$\bullet$ the Gauss map may be identified with
$\pi|_S$, and hence has degree $\ell$;

$\bullet$ the total Gaussian
curvature $\int\!\!\int {\cal K} ds^2 $ of the induced branched metric equals
$-4\pi\ell$;

$\bullet$ if $S$ is $\tau$-invariant then $\phi$ factors through
$S/\tau$; 

$\bullet$  if $G\subset {\rm SO(3)}$ is the symmetry group of a
regular solid and $\tilde{G}\subset{\rm SU(2)}$ the  corresponding binary
group  and $S$ is invariant under $\tilde{G}$ then $G$ is a subgroup of the
symmetry group  of the corresponding minimal surface.     

{\it Remark.}
$\Omega$ induces the branched  metric $4ds^2$ on $S^\ast$.   

The geometric correspondence described above underlies  the
{\it Weierstrass formulae in free form}.   In `global form', in
which $S$ is described by a pair of meromorphic functions   on a Riemann
surface: $(g, f):M\ar\C^2$, the coordinate functions of the null curve  
$\Om:M^\ast\ar\C^3$  are  given by: 
 \begin{eqnarray} 
 \Om_1 &=&
\frac{1}{2}\left(-\frac{1}{2}(1-g^2)\frac{d^2f}{dg^2}-g\frac{df}{dg}+f\right)\\[0.1in]   \Om_2 &=&
\frac{i}{2}\left(-\frac{1}{2}(1+g^2)\frac{d^2f}{dg^2}+g\frac{df}{dg}-f\right)\\[0.1in] 
 \Om_3 &=& \frac{1}{2}\left(g\frac{d^2f}{dg^2}-\frac{df}{dg}\right) 
\end{eqnarray} 
 where $\di{\frac{df}{dg}=\frac{f'}{g'}}$ and
$\di{\frac{d^2f}{dg^2}=\left(\frac{df}{dg}\right)'\frac{1}{g'}}$, etc.     
$\phi=$Re($\Om$)  describes a branched minimal   surface in $\R^3$.   Note 
that $g$  may be  identified with the classical Gauss map of $\phi$.

{\it Remark.} These formulae are   not  canonical: their precise shape is
determined by the real structure  $\tau$, and thus the choice of  $\R^3$ in
$\gs$ made in \cite{h1}.   They  differ slightly from the classical formulae,
\cite{da}. (Note however that the formulae   given in the appendix of 
\cite{h1} require slight adjustment.)    

 The following is an immediate
consequence of the nature of osculation duality. 

{\bf Proposition 2.1} {\it
Suppose that $S\subset\T$ is the spectral curve of a monopole and
$\Omega:S^*\ar\C^3$, the associated null curve. Then the affine null planes in
$\C^3$ that osculate the null curve intersect $\R^3$ in the spectral lines of
the monopole.} 

{\it Remark.} This requires clarification for charge 1,
since in that case the spectral curve is a section, and osculation of it is
degenerate, in the sense that it gives  only a point in $\C^3$. The affine
null planes through a point in $\C^3$ should be viewed as `osculating that
point'. Similar remarks apply to any spectral curve which includes a global
section as a component.


\vspace*{0.3in}

{\bf \S 3. Osculation of spectral curves of charge 2 monopoles.} 
 
In this section we describe  the main features of the minimal surfaces 
generated by osculation of charge 2 monopole spectral curves.  At this point we refrain
from deriving formulae for the minimal surfaces: these are discussed in the
next section.    

  Hurtubise \cite{hurt2}, showed that the spectral curve $S$, of a centred
charge   2 monopole, may  by rotation of $\R^3$, be brought to the reduced
form:  
\begin{eqnarray} 
 \eta^2=r_1\z-r_2\z^2-r_1\z^3, \ \ r_1,\ r_2\in \R,\
r_1\geq 0.  
\end{eqnarray} 
When $r_1\neq 0$, the triviality of $\rl^2|_S$ constrains the real period: $\om_1=2\sqrt{r_1}$.  
At $r_1=0$, $S$ degenerates into the pair of global sections given by 
$\eta^2=-\pi^2\z^2/4$: this gives the reduced form for an axially symmetric
monopole.    Observe that if $r_1\neq 0$, then $S$ is a smooth 
$\tau$-invariant elliptic curve on $\T$.  $\tau$ restricted to $S$ takes the
form $\tau(u)=-\overline{u}+\om_3/2$, and has no fixed points: the associated
lattice is rectangular.      

When $r_1\neq 0$, the branch points of $S$ are at $0,\ \infty,\ -a,\ {\rm and} \ a^{-1}$,
 where $-a$ and 
 $a^{-1}$ are the roots of  $\z^2-\frac{r_2}{r_1}\z-1=0$.  These two antipodal pairs give 
 $\alpha$  and $\beta$, the spectral lines through 0, the centre of the monopole. 
 The monopole has a distinguished  bisector $\un{e}_1$, of the angle between
$\alpha$ and $\beta$: this is the monopole's {\it main axis}. (In the axially
symmetric case this is the axis of symmetry.) The second bisector, $\un{e}_2$,
 is the monopole's {\it Higgs axis}. The perpendicular through 0 to $\un{e}_1$
and $\un{e}_2$, is denoted $\un{e}_3$, and called the {\it third axis}. Having
fixed $\un{e}_1$, $\un{e}_2$ and $\un{e}_3$,  the monopole is determined by an
angle $0\leq\theta<\pi/2$. ($\theta=0$, gives a centred axially symmetric
monopole.) See \cite{ah} for further details.

It is observed in
\cite{hurt2} that $S$  has symmetries,  permuting the roots of (4): these
correspond to the subgroup $D$ of SO(3), comprising rotations through $\pi$
about the axes $\un{e}_1$, $\un{e}_2$ and $\un{e}_3$ in $\R^3$. $D$ is of
course isomorphic to ${\Bbb Z}_2\times{\Bbb Z}_2$.

Following \cite{ah},  let   $\tan(2\theta)=\frac{2r_1}{r_2}$. $2\theta$ is the
angle between the lines $\alpha$ and $\beta$.  It is natural and eases
calculation to introduce the modulus $k=\sin(\theta)$, together with the
complementary modulus $k'=\cos(\Th)$.     SO(3) acts naturally on the moduli
space ${\cal M}_2^0$ of centred 2-monopoles. The orbits are parameterized by
$\theta$, or equivalently, $k\in[0,1)$. For $\theta=0$, the orbit is
isomorphic to $\R\p_2$: this parameterizes  the centred axially symmetric
2-monopoles. For $\theta\neq0$, the orbit is isomorphic to  SO(3)/$D$.     In
\cite{ah}, it is observed that the triviality of $\rl^2|_S$ means that (4) may
be rewritten:  \begin{eqnarray}  \eta^2=K(k)^2\z(kk'(\z^2-1)+(k^2-(k')^2)\z), 
\end{eqnarray}  where as usual,
$\di{K(k)=\int_0^{\pi/2}\frac{d\psi}{\sqrt{1-k^2\sin^2\psi}}}$. Accordingly,
emphasising $k$ dependence, from now on  we refer to this curve as $S_k$.

    Now, a global
section $\sigma_z$ of $\T$, corresponding to $z\in\C^3$, osculates $S_k$ at $p$
if and only if $\sigma_{\overline{z}}$ osculates $S_k$ at $\tau(p)$, thus
$\Omega(\tau(u))=\overline{\Omega(u)}$, and hence $\phi(\tau(u))=\phi(u)$. A
cursory inspection of the structure of the Weierstrass formulae for these
surfaces shows that $\phi(-u)=-\phi(u)$, (cf. Theorem 4.7). These mean that
$\phi$ enjoys many  symmetries on the fundamental period rectangle. In
particular observe that $\phi$ is defined on the doubly punctured Klein bottle
$S_k/\tau-\{[0], [\om_3/2]\}$.  

  As a curve on ${\cal C}(Q)$ in
$\p_3$, $S_k$ has degree 4. The points of hyperosculation on $S_k$ are the points
of order 4 in its group structure: this follows from Abel's theorem. Each of
the four branch points of $\pi|_{S_k}$, is a point of hyperosculation since the
osculating hyperplane at a branch point $b$ say, lies tangent to ${\cal C}(Q)$
along the fibre through  $b$,  and thus  intersects $S_k$  with multiplicity 4
at $b$. These give the points of order 2 in the group structure of $S_k$. They
come in two antipodal pairs and correspond to the  two ends of the minimal
surface  in $\R^3$. This leaves twelve zeros in the branched metric on
$S_k^\ast$: these pass to six branch points on $S_k/\tau-\{[0], [\om_3/2]\}$. 
   Comparing these observations with the properties of osculation duality
listed in $\S2$ gives:  

{\bf Proposition 3.1} {\it (i) Osculation of the spectral
curve $S_k$ of a non-axially symmetric centred 2-monopole gives a branched
minimal immersion of the punctured Klein bottle $\phi:S_k/\tau-\{[0],
[\om_3/2]\}\ar\R^3$, with the following properties: 

(ii) the total
Gaussian curvature of the induced  branched metric on $S_k^\ast$ equals
$-8\pi$. 
 
(iii) The minimal surface has two ends. These are perpendicular to the two
spectral lines through the monopole's centre. 

(iv) There are six branch
points, (of ramification index 1), on the minimal surface in $\R^3$. These
are:   $$  \pm \beta_1=\pm\phi(\om_1/4),\  \pm\beta_2=\pm\phi(\om_2/4)\ and\ 
\pm\beta_3=\pm\phi(\om_3/4).  $$ 

(v)  The image of $\phi:S_k/\tau-\{[0],
[\om_3/2]\}\ar\R^3$,  is invariant under $D$.} 

{\it  Remark.} Osculation
of the (reduced and centred) axially symmetric 2-monopole spectral curve
yields the pair of points $(0,0, \pm \frac{i\pi}{4})$ in $\C^3$. (So the
auxiliary `minimal surface' in this case is the point at the origin.)

Locally, around each of the points of hyperosculation that are not branch
points of $\pi|_{S_k}$, $S_k$ may be described by $\eta=a_4\z^4+{\cal
O}(\z^5)$, for some $a_4\in\C$. Hence at each of the branch points the minimal
surface is locally a perturbation of a rescaled associate surface of the
minimal surface determined by osculation of $\eta=\z^4$.     It is easy to see
directly from calculation  that the latter maps the lines $x=0$,
$y=\frac{1}{\sqrt{3}}x$, and $y=-\frac{1}{\sqrt{3}}x$,  to the three rays in
the $(x_1,x_2,0)$-plane at $120^o$, where $x=0$, is mapped $2:1$ onto 
$\{(x_1,0)\ ;\ \ x_1\leq0\}$, etc. Each of the $60^o$ sectors in the
$(x,y)$-plane is embedded, along with the image of its reflection through the
origin, to a surface bounded by two of the rays, cf. Figure 6 in
\cite{lawson}.  This {\it triple curve intersection structure} at the branch
point is stable under higher order perturbations and  multiplication by $a_4$.
(Of course, it `twists' if $a_4\not\in\R$.) Hence we see six of these triple
curve intersection structures at the branch points on the monopole minimal
surface.

\newpage

{\bf \S 4.  Formulae for the null curve.}

 In \cite{hurt2}, the following substitutions  are introduced for  $r_1\neq
0$:  
$$
\z=\tilde{\z}+k_2 \ \ {\rm and}\ \  \eta=k_1\tilde{\eta},
$$
  where 
$k_1=\frac{1}{2}\sqrt{r_1}$ and $k_2=r_2/3r_1$. Thus (4) becomes 
\begin{eqnarray}
  \tilde{\eta}^2=4\tilde{\z}^3-g_2\tilde{\z}-g_3, 
\end{eqnarray}
  where $g_2=12k_2^2+4$ and $g_3=8k_2^3+4k_2$.    

If $\wp(u)$ is the Weierstrass $\wp$-function determined by $g_2$ and $g_3$,
then the spectral curve    $S$ is uniformised by   $\z=g(u)=\wp(u)+k_2$ and
$\eta=f(u)=\frac{\om_1}{4}\wp'(u)$. 

{\it Remark.}  It should be noted that
direct substitution of these  into the Weierstrass formulae (1)-(3), yields
very complicated expressions. We now outline another approach which results in
the relatively simple formulae described in Theorem 4.7 below.    First the
meaning of the parameter $k_2$ is  clarified: 
 
{\bf Proposition 4.1} $k_2=-e_3$.  

{\it Proof.}  $\wp'(\om_j/2)=0$, implies
$4e_j^3-g_2e_j-g_3=0$,  for $j=1,2,3$, and hence  $$ 
(e_j+k_2)(e_j^2-k_2e_j-(1+2k_2^2))=0,\ \ {\rm for}\ j=1,2,3.  $$  The roots of
the quadratic factor are $(k_2\pm\sqrt{9k_2^2+4}\ )/2$. The ordering of the
roots   follows from the elementary fact   that for a rectangular lattice
$\wp(u)$   takes real values, and is strictly decreasing as $u$ passes around
the rectangle with vertices   $0,\ \om_1/2,\ \om_3/2,\ \om_2/2$, and hence
$e_1>e_3>e_2$, see \cite{duval}. It follows that $k_2=-e_3$. $\ \Box$  

  {\bf Corollary 4.2} {\it  The half-period values of $\wp(u)$ are given by:} 
\begin{eqnarray}
  e_1=\frac{2-k^2}{3kk'},\ \ e_2=-\frac{1+k^2}{3kk'} \ \ {\rm
and} \ \ e_3=\frac{2k^2-1}{3kk'}. 
 \end{eqnarray}    

{\it Proof.} From the definitions of $\theta$ and $k_2$, it follows   that
$\di{k_2=\frac{1-2k^2}{3kk'}}$. The result is immediate.  $\ \Box$ 

{\it Remarks.} (i)  A simple calculation shows that   
\begin{eqnarray}  g_2=
\frac{4(1-k^2+k^4)}{3k^2k'^2}\ \ \ {\rm and}\ \
g_3=\frac{4(k^2-2)(k^2+1)(2k^2-1)}{27k^3k'^3}.  
\end{eqnarray}    

(ii) The periods of $S_k$ are given
by:  \begin{eqnarray}  \om_1=2\sqrt{kk'}K(k)\ \ {\rm and}\ \
\om_2=2i\sqrt{kk'}K'(k),  \end{eqnarray}  where,\ as usual, $K'(k)=K(k')$,
\cite{ah}.   

 {\bf Lemma 4.3}\ {\it For $\wp(u)$, with $g_2,\ g_3$ determined
by $k$ as above:}  \begin{eqnarray} 
1+2\di{\frac{k'}{k}}(\wp(u)-e_3)-(\wp(u)-e_3)^2&=&\wp'(u)\sqrt{\wp(2u)-e_1}\\[0.15in]  -1+2\di{\frac{k}{k'}}(\wp(u)-e_3)+(\wp(u)-e_3)^2&=&-\wp'(u)\sqrt{\wp(2u)-e_2}\\[0.15in]   1+(\wp(u)-e_3)^2&=&-\wp'(u)\sqrt{\wp(2u)-e_3}  \end{eqnarray} 
{\it where, in the first formula, $i$ times the radical is positive at
$u=\om_2/4$, and in the other two, the radical is positive at $u=\om_1/4$.}   
 
{\it Proof.} The quarter-period formulae  of the Appendix, together with the translation formulae,
 yield: 
$$ 
\begin{array}{lclcl} 
\di{\wp(\frac{\om_1}{4})-e_3}&=&
\di{\wp(\frac{3\om_1}{4})-e_3}&=&\di{{\frac{1+k'}{k}},} \\[0.15in] 
\di{\wp(\frac{\om_1}{4}+\frac{\om_2}{2})-e_3}&=&
\di{\wp(\frac{3\om_1}{4}+\frac{\om_2}{2})-e_3}&=&\di{{\frac{-1+k'}{k}}.}
\end{array}  $$  Now observe that  $(1+k')/k$, $(-1+k')/k$ are the roots of  
$$ 
1+2\frac{k'}{k}x-x^2=0. 
$$ 
 It is  clear that both sides in (10) have poles of order $-4$ at 0, and zeros of order 1  
at $\om_1/4$, $3\om_1/4$, $\om_1/4+\om_2/2$ and $3\om_1/4+\om_2/2$. 
The scaling is fixed  at $\om_2/4$, cf. the Appendix. (11) and (12) are proved
similiarly.\ $\Box$   

{\bf Definition.}  Let $\ff_j(u)$ denote the square root of $\wp(u)-e_j$, whose residue 
at the origin is $1$. 
 
These choices of signs accord with those in Lemma 4.3. 
 
{\it Remark.}  This notation follows \cite{duval}.

{\bf Lemma 4.4}\ \ {\it Substitution  of}\ \ $g(u)=\wp(u)-e_3$,
$\di{f(u)=\frac{\om_1}{4}\wp'(u)}$,  {\it into (2) yields:} 
$$ 
\Omega_2(u)=\frac{i\om_1}{4}\left\{\frac{1+(\wp(u)-e_3)^2}{\wp'(u)}\right\}^3. 
$$ 
 
{\it Proof.}  First observe that  
\begin{eqnarray*} 
g\frac{df}{dg}-f&=&\frac{\om_1}{4}\left\{\frac{gg''-(g')^2}{g'}\right\}\\[0.1in] 
                &=&\frac{\om_1}{4\wp'}\left\{2\wp^3-6e_3\wp^2+
\frac{g_2}{2}\wp+g_3+e_3\frac{g_2}{2}\right\}\\[0.1in] 
                &=&\frac{\om_1(1+g^2)g}{2g'}, 
\end{eqnarray*} 
since $g_2=4(1+3e_3^2)$, and $g_3=-4e_3(1+2e_3^2)$.  
 
Now observe that   $\di{\frac{d^2f}{dg^2} = \frac{\om_1}{4}\left\{\frac{g'''g'-(g'')^2}{(g')^3}\right\}}$, and hence 
\begin{eqnarray*} 
            \Omega_2 &=&
\frac{i\om_1(1+g^2)}{16g'^3}\left\{4g(g')^2-g'''g'+(g'')^2\right\}\\[0.1in]                       &=&
\frac{i\om_1(1+g^2)}{16g'^3}\left\{4\wp^4-16e_3\wp^3+2g_2\wp^2+4(2g_3+g_2e_3)\wp+4e_3g_3+\frac{g_2^2}{4}\right\}\\[0.1in]               
       &=& \frac{i\om_1}{4}\left\{\frac{1+g^2}{g'}\right\}^3,  \end{eqnarray*} 
again using $g_2=4(1+3e_3^2)$, and $g_3=-4e_3(1+2e_3^2)$. $\ \Box$

{\bf Lemma 4.5}  $$\di{\Omega_2(u)=-i\frac{\om_1}{4}\ff_3(2u)^3}.$$ 
 
{\it Proof.}\ This follows immediately from (12) and  4.4.\ $\Box$  
 
{\bf Lemma 4.6}\  {\it For $g(u)=\wp(u)-e_3$, and
$\di{f(u)=\frac{\om_1}{4}\wp'(u)}$:}  \begin{eqnarray} 
\frac{d^3f}{dg^3}(u)=-3\om_1\frac{\wp'(2u)}{\wp'(u)^2}. 
\end{eqnarray} 
{\it Proof.} Differentiate (2) and the formula of Lemma 4.5 and compare,
using (12).$\ \Box$     
 
With respect to the {\it $k$-monopole coordinates} $(\un{e}_1, \un{e}_2,
\un{e}_3)$, the null curve is represented by $\Phi=A_k\Omega$, i.e. 
\begin{eqnarray} 
 \left(\begin{array}{c}\Phi_1\\ \Phi_2\\
\Phi_3\end{array}\right)=\left(\begin{array}{rcl}-k & 0 & k'\\ k'&0&k\\
0&1&0\end{array}\right)\left(\begin{array}{c}\Om_1\\ \Om_2\\
\Om_3\end{array}\right) 
 \end{eqnarray} 
\newpage
{\bf Theorem 4.7} {\it Suppose that
$k\in(0,1)$, determines $g_2,\ g_3$, as in (8), and $\wp(u)$ is the associated
Weierstrass function. The null curve that is  generated by osculation of the 
spectral curve   described implicitly by $g(u)=\wp(u)-e_3$,
$\di{f(u)=\frac{\om_1}{4}\wp'(u)}$,   has components with respect to 
$k$-monopole coordinates given by:}  \begin{eqnarray}  \Phi_1(u) &=&  
-k\frac{\om_1}{4}\ff_1(2u)^3 \\[0.1in]  \Phi_2(u) &=&  \ \
k'\frac{\om_1}{4}\ff_2(2u)^3  \\[0.1in]  \Phi_3(u) &=&
-i\frac{\om_1}{4}\ff_3(2u)^3  \end{eqnarray}     
{\it  Proof.} \ The formula for $\Phi_3(u)$ is equivalent to that of Lemma
4.5.    
Differentiating with respect to $u$: 
\begin{eqnarray*} 
\Phi_1'(u)&=&g'(u)\frac{d\Phi_1}{dg}(u)\\[0.15in] 
          &=&\frac{k}{4}g'(u)\frac{d^3f}{dg^3}(u)(1-g(u)^2+2\frac{k'}{k}g(u))\\[0.15in] 
          &=&-\frac{3}{4}k\om_1\wp'(2u)\ff_1(2u) 
\end{eqnarray*} 
where the last equality follows from (10), together with Lemma 4.6. Therefore 
$$\Phi_1(u)=-k\frac{\om_1}{4}\ff_1(2u)^3+\Phi_1(\frac{\om_1}{4}).$$ 
Now note that it follows from the quarter-period formulae (of the Appendix) that $\Phi_1(\om_1/4)=0$. 
 
Similiarly, 
\begin{eqnarray*} 
\Phi_2'(u)&=&g'(u)\frac{d\Phi_2}{dg}(u)\\[0.15in] 
          &=&\frac{k'}{4}g'(u)\frac{d^3f}{dg^3}(u)(-1+g(u)^2+2\frac{k}{k'}g(u))\\[0.15in] 
          &=&\frac{3}{4}k'\om_1\wp'(2u)\ff_2(2u) 
\end{eqnarray*} 
where the last equality follows from (11), together with Lemma 4.6. Therefore 
$$\Phi_2(u)=k'\frac{\om_1}{4}\ff_2(2u)^3+\Phi_2(\frac{\om_2}{4}).$$ 
Finally note that it follows from the quarter-period formulae (of the Appendix) that
 $\Phi_2(\om_2/4)=0$.\ $\Box$ 
 
{\it Remark.}  The nullity of $\Phi$ is a special case of the quadratic identity 
$$ 
(e_2-e_3)\ff_1^2(z)+(e_3-e_1)\ff_2^2(z)+(e_1-e_2)\ff_3^2(z)=0. 
$$ 
For a normal rectangular lattice this reduces to the following identity between Jacobi 
functions: $-k^2cs^2(z)-k'^2ns^2(z)+ds^2(z)=0$. Cf. $\S 4$ of \cite{duval}. 
 
{\bf Corollary 4.8}\ {\it The  branched metric induced on the spectral curve
by $\phi$ has the form:}  $$ 
ds^2(u)=\frac{9}{32}|\wp'(2u)|^2\{k^2|\wp(2u)-e_1|+k'^2|\wp(2u)-e_2|+|\wp(2u)-e_3|\}|du|^2. 
$$ 
 
{\it Remark.} Observe that this displays the branching and end structure at the quarter-period
 points  described in $\S 3$. 
 

\vspace{0.3in}

{\bf \S 5. The locations of the branch points.} 
 
The following can be derived from Theorem 4.7, (or directly from
quarter-period formulae). They are written w.r.t. $k$-monopole coordinates.
  
{\bf Proposition 5.1} 
\begin{eqnarray} 
\Phi(\frac{\om_1}{4})&=&\frac{K(k)}{2k}(0,\ 1,\ -ik'^2)\\[0.1in] 
\Phi(\frac{\om_2}{4})&=&\frac{K(k)}{2k'}(-i,\ 0,\ k^2)\\[0.1in] 
\Phi(\frac{\om_3}{4})&=&\frac{K(k)}{2}(-ik'^2,\ k^2,\ 0) 
\end{eqnarray} 
{\bf Corollary 5.2} {\it Projecting to $\R^3$ gives:} 
 \begin{eqnarray} 
\be_1&=&\frac{K(k)}{2k}(0,\ 1,\ 0)\\[0.1in] 
\be_2&=&\frac{k^2K(k)}{2k'}(0,\  0,\ 1)\\[0.1in] 
\be_3&=&\frac{k^2K(k)}{2}(0,\ 1,\ 0) 
\end{eqnarray} 
{\it Remark.} Observe that the positions of the branch points are tied to the moduli space 
parameter $k$ as above, through the non-singularity constraint ${\rm
L}^2|_{S_k}={\cal O}_{S_k}$.   
 There are a number of interesting consequences which we list in the following corollary. First note the following 
 
{\bf Definition.}\  For values of $k$ close to 1, by  {\it separation distance} we mean $K(k)$.

{\bf Corollary 5.3} {\it (a) The branch points $\pm\be_1$, $\pm\be_3$,  lie
on the Higgs axis, $\underline{e}_2$.
   
(b) In the `widely separated limit', $k\ar 1$: 
 
$\bullet$ $\pm\be_1,\ \pm\be_3 \ar\pm\frac{K(k)}{2}(0,\  1,\  0) $, respectively.

Thus, these pairs of branch points approach  the `star centres' exponentially fast relative to the separation distance.  
 
$\bullet$ $\pm\be_2\ar \pm\infty$  along $\underline{e}_3$, exponentially fast relative to separation distance. 
 
(c) As $k\ar 0$, and  the monopoles approach the axially symmetric solution: 
 
$\bullet$ $\pm\be_2,\ \pm\be_3\ar 0$. 
 
$\bullet$ As $k\ar0$, $\pm\be_1$ approach the origin (with $\pm\be_3$, respectively,) but `turn around', at $k_0$, such that $k_0dK/dk(k_0)=K(k_0)$, and go back out to $\pm\infty$, respectively,  along the Higgs axis. } 
 
{\it Remark.}\ $k_0<1/\sqrt{2}$, the value that gives the square lattice. 
 
{\bf Definition.} We call $\pm\be_1$, $\pm\be_3$,   the {\it Higgs branch points}. 
 
{\it Remark.}\  It is instructive to take note of the asymptotic behaviour of the imaginary components at the branch points.

 
\vspace{0.3in}

 {\bf \S 6. Curves of symmetry of $\phi={\rm Re}(\Phi)$.} 
 
{\bf Definition.} 
 
(i) $V_a:= \{a\om_1+iy\ ;\ 0\leq y<|\om_2|\}$. 
  
(ii) $H_a:= \{x+a\om_2\ ;\ 0\leq x<|\om_1|\}$. 
 
We study  $\Gamma_{S_k}=\phi(H_{\frac{1}{4}}\cup V_{\frac{1}{4}}\cup
V_{\frac{3}{4}})$: this is a graph on the minimal surface 
which connects the six  branch points.

{\bf Definition.} Write
$V_{\frac{1}{4}}=V_1\cup V_2\cup V_3\cup V_4$ where
$V_j:=\{\frac{\om_1}{4}+iy\ ;\ (j-1)\frac{|\om_2|}{4}\leq y\leq
j\frac{|\om_2|}{4}\}$,  $j=1,2,3,4$.

{\bf Proposition 6.1} \ {\it  The
branched minimal immersion  $\phi$ maps $V_{\frac{1}{4}}\cup V_{\frac{3}{4}}$
to  the Higgs axis so that:    (i) $\phi|_{V_3}=\phi|_{V_1}$ and
$\phi|_{V_4}=\phi|_{V_2}$.   
(ii) $\phi$ maps $V_1$ monotonically from $\beta_1$ to $\beta_3$. 
 
(iii) $\phi$ maps $V_2$ monotonically from $\beta_3$ to $\beta_1$. 
 
(iv) $\phi(\frac{3\omega_1}{4}+iy)=-\phi(\frac{\omega_1}{4}+i(\frac{\omega_2}{2}-y))$.} 
 
{\it Proof.} \ (i) follows immediately from $\phi(\tau(u))=\phi(u)$. 
 
(ii) and (iii): for $u\in V_1\cup V_2$, $e_3 \leq\wp(2u)\leq e_1$. Hence, for such $u$,  $\phi_1(u)\equiv 0$, and $\phi_3(u)\equiv 0$. $\phi_2(u)$ monotonically decreases (through real values) from $k'(e_1-e_2)^{3/2}$ to $k'(e_3-e_2)^{3/2}$ between $\om_1/4$ and $\om_3/4$. It is monotonically increasing between $\om_3/4$ and $\om_1/4+\om_2/2$. 
 
(iv) follows immediately from $\phi(-\tau(u))=-\phi(u)$. \ $\Box$ 
 
{\bf Proposition 6.2} \ {\it $\phi$ maps $H_{\frac{1}{4}}$ monotonically onto
 the `star' in the $(\un{e}_2, \un{e}_3)$-plane, with vertices $\pm \beta_2,\
\pm \beta_3$, that is formed when $D$ acts on the concave curve \newline
$\phi(\{x+\om_2/4\ ;\ 0\leq x\leq \om_1/4\})$.}  

  {\it Proof.}\ Since
$\wp(2u)-e_1<0$, for all $u\in H_\frac{1}{4}$, it follows from 4.7 that
$\phi_1(u)\equiv 0$ on $H_\frac{1}{4}$.    $\phi_2$ is odd at $\om_2/4$ and
$\om_2/4+\om_1/2$, while it is even at $\om_3/4$ and $\om_3/4+\om_1/2$. On the
other hand, $\phi_3$ is odd at $\om_3/4$ and $\om_3/4+\om_1/2$, while it is
even at $\om_2/4$ and $\om_2/4+\om_1/2$. This leads immediately to
$D$-invariance.    For $u\in H_\frac{1}{4}$, $\wp'(2u)$,
$(\wp(2u)-e_2)^{\frac{1}{2}}$ and $(e_3-\wp(2u))^\frac{1}{2}$ are real. Hence,
for $u\in\{x+\om_2/4\ ;\ 0\leq x\leq \om_1/4\}$,  $$ 
\frac{d\phi_3}{d\phi_2}(u)=\frac{{\rm Re}(\Phi_3'(u))}{{\rm
Re}(\Phi_2'(u))}=-\frac{1}{k'}\sqrt{\frac{e_3-\wp(2u)}{\wp(2u)-e_2}}.  $$  For
such $u$, $e_3-\wp(2u)$ monotonically decreases from $e_3-e_2$ to $0$,   while
$\wp(2u)-e_2$ monotonically increases from $0$ to $e_3-e_2$. Thus observe that
the derivative  monotonically increases from $-\infty$ to $0$ on the curve 
between $\beta_2$ and $\beta_3$.\ $\Box$  

  {\bf Definition.} Write
$\Gamma_{S_k}=\Gamma_{\rm Star}(k)\cup\Gamma_{\rm Higgs}(k)$, where $\Gamma_{\rm
Star}(k)=\phi(H_\frac{1}{4})$.

    Note that $\Gamma_{S_k}$ is $D$-invariant.   
Similiar arguments give:   

 {\bf Proposition 6.3}\ (i) $\phi$ maps $V_0$ onto
the $\un{e}_3$-axis between $\beta_2$ and $\infty$. Similiarly, $\phi$ maps
$V_{\frac{1}{2}}$ onto the $\un{e}_3$-axis between $-\beta_2$ and $-\infty$.

  (ii) Near $\beta_2$, $\phi(V_0)$, together with the curves of $\Gamma_{\rm
Star}$ that emanate from $\beta_2$, comprise the triple curve intersection
structure at $\beta_2$. The analogue holds at $-\beta_2$.  

   (iii) The triple
curve intersections  at $\pm\beta_3$ are induced by $\Gamma_{S_k}$.  

  (iv) $\phi$
maps $H_0$ between $0$ and $\om_1/2$, to two arms of the triple curve
intersections at $\beta_1$. (The third arm is $\phi(V_\frac{1}{4})$.) These
two arms emanating from $\beta_1$ lie in the $(\un{e}_1, \un{e}_2)$-plane and
are concave down. The analogue holds at $-\beta_1$.


 \vspace{0.3in}

{\bf \S7. Spectral lines near $\Gamma_{\rm Higgs}(k)$, as $k\ar1$.}

 Asymptotically, as $k\ar1$, $K(k)\sim -{\rm log}k'$. Thus
observe from (5), that as $k\ar1$, the configuration of spectral lines
approximates the two stars on the Higgs axis through the points at distance
$K(k)/2$ from the origin, see \cite{ah}. In 5.2, we saw that $\Gamma_{\rm Higgs}(k)$ 
shrinks to these points as $k\ar1$. Here we observe that in this limit the
stars approximate the normal lines to the minimal surface in the vicinity of
$\Gamma_{\rm Higgs}(k)$.

The branch points on the null curves in $\C^3$ which project to the Higgs
branch points approach exponentially close to $\R^3$, relative to separation
distance, as $k\ar1$. In fact the same is true  at all points along
$\Gamma_{\rm Higgs}(k)$, in this limit: this can be deduced from 4.7. In the
light of 2.1, elementary arguments yield:

{\bf Theorem 7.1} {\it As $k\ar1$, every normal line to the minimal surface,
along $\Gamma_{\rm Higgs}(k)$, becomes exponentially close, relative to
separation distance,  to a spectral line of the monopole.}

{\it Remark.}  Of course, this approximation is true anywhere `close enough'
to $\Gamma_{\rm Higgs}(k)$. In general, far away from $\Gamma_{\rm Higgs}(k)$, the
normal lines are not close to spectral lines. Consider, for example, the fact
that the  maximum distance from the origin, for fixed $k$, of any spectral
line is finite. However, the minimal surface has flat ends and out on these
will be normal lines at arbitrary distance from the origin.     

 Note in
particular that the latter observation applies at the two branch points
$\pm\beta_2$, on the monopole's third axis. This is because the branch points
on the null curve in $\C^3$ which engender these move away from $\R^3$, as
$k\ar1$.


\vspace{0.3in}  
 
{\bf \S 8. The Gauss map and curvature concentration.} 
 
Recall that the spectral curve of a  charge $\ell$ monopole is an $\ell$-fold
branched cover of $\p_1$ in $\T$. Since the projection map to $\p_1$ may be
identified with the Gauss map, $g$, of the auxiliary minimal surface
determined by osculation, it  follows that the monopole energy is:  $$  {\cal
E}(\nabla, \Phi)=4\pi{\rm deg}(g)=4\pi \ell.  $$   Our purpose here is to 
interpret  this global coincidence `locally',  when $\ell=2$.    
Let ${\cal K}$ denote the Gaussian curvature of the branched metric $ds^2$ on $\C/(\Z\om_1+\Z\om_2)$ induced by the branched minimal immersion. Recall that  
$$ 
\int_{\C/(\Z\om_1+\Z\om_2) }{\cal K}ds^2=-\int_{\C/(\Z\om_1+\Z\om_2) }\frac{4|g'|^2}{(1+|g|^2)^2}dxdy=-4\pi \ell, 
$$ 
cf. \cite{lawson}. The second integral is just the area induced by the Gauss map $g$: we 
study the behaviour of  
$$ 
G=\di{\frac{4|g'|^2}{(1+|g|^2)^2}}, 
$$ 
in the limits $k\ar0,\ 1$.  In particular, we show that it {\it localises} in these limits.
 This reflects the behaviour of ${\cal K}$  on the surface in $\R^3$: as
$k\ar 1$, ${\cal K}$ localises at the monopole particles. It follows from
Theorem 7.1 that in the limit $k\ar1$, ${\cal K}$ measures the twisting of the
spectral lines through the particles. The total twisting, measured in the
induced metric, equals the monopole's energy.   

 Recall that the Gauss map is
given by $g(u)=\wp(u)-e_3$. This is related to the Euclidean Gauss map
$\gamma$, of the minimal surface, via stereographic projection:
$\gamma=P^{-1}\circ g$, where $P:S^2\ar\C\cup\{\infty\}$,
$P(x)=(x_1+ix_2)/(1+x_3)$. $g$ maps the rectangle with vertices $0, \om_1/2,
\om_3/2, \om_2/2$, to the lower half plane so that its boundary maps to $\R$.
Obviously, the behaviour of $g$ on the whole period domain may be inferred
from this. The  behaviour of the two pairs:  $$\begin{array}{ll}  g(0)=\infty,
& g(\om_3/2)=0,\ {\rm and}\\[0.1in]  g(\om_1/2)= k'/k, & g(\om_2/2)=-k/k', 
\end{array}$$  corresponding to the two spectral lines through the origin,
elucidates the manner in which $g$ covers $\p_1$. In particular, as $k\ar1$,
it appears that $G$ `concentrates' on the quarter-period lines $x=\om_1/4$,
and $x=3\om_1/4$.  The lines $x=0$, and $x=\om_1/2$, map to the circles on
$S_k$ that are `asymptotically pinched off', as $k\ar1$. Moreover, as
$k\ar0$, $G$ appears to `concentrate' on the quarter-period lines
$iy=\om_2/4$, and $iy=3\om_2/4$. The lines $iy=0$, and $iy=\om_2/2$, map to
the circles on $S_k$ that are pinched off at $k=0$.     Despite the simplicity
of the expression for $g$,  a deeper insight into these matters is gained
through working in monopole coordinates: the symmetries of $G$ with respect to
the quarter-period lines, and its behaviour on them are thus revealed.         

 The Gauss map
$\gamma_\Phi:\C/(\Z\om_1+\Z\om_2)\ar Q_1$, with respect to monopole coordinates, is given by
differentiating (15)-(17):  $$  \gamma_\Phi(u)=[-k\ff_1(2u),\ k'\ff_2(2u),\
-i\ff_3(2u)].  $$  It follows that $g_\Phi:\C/(\Z\om_1+\Z\om_2)\ar\C\cup\{\infty\}$, is given
by:  $$ 
g_\Phi(u)=h^{-1}\circ\gamma_\Phi(u)=\frac{-i\ff_3(2u)}{k\ff_1(2u)+ik'\ff_2(2u)},  $$ 
where $h:\C\cup\{\infty\}\ar Q_1$, $h(\z)=[1-\z^2, i(1+\z^2), -2\z]$. 

  {\it Remark.} Let $\gamma_\phi:\C/(\Z\om_1+\Z\om_2)\ar S^2$, be the Euclidean
Gauss map of the minimal surface $\phi={\rm Re}(\Phi)$. Observe that $g_\Phi$
agrees with $\gamma_\phi$,  after the latter is composed with stereographic
projection from $-\un{e}_3$ to the $(\un{e}_1, \un{e}_2)$-plane in $\R^3$.    

  {\bf Proposition 8.1}{\it (i)  $$ 
g_\Phi(u)=-\frac{k'\ff_2(2u)+ik\ff_1(2u)}{\ff_3(2u)}.  $$ 
(ii)} 
$$g_\Phi'(u)=\frac{2ig_\Phi(u)}{\ff_3(2u)}. 
$$ 
{\it Proof.} (i) follows immediately from $k^2\ff_1^2+k'^2\ff_2^2-\ff_3^2=0$. 
 
(ii) 
$$ 
g_\Phi'(u)=-(\frac{\wp'}{\ff_3^3}(\ff_3^2(k'\ff_2^{-1}+ik\ff_1^{-1})-(k'\ff_2+ik\ff_1)))(2u), 
$$ 
and hence $\ff_3^2=(k'\ff_2+ik\ff_1)(k'\ff_2-ik\ff_1)$, implies 
\begin{eqnarray*} 
g_\Phi'(u)&=&(\frac{\wp'}{\ff_3^2}(k'^2+k^2-1+ikk'(\ff_1^{-1}\ff_2-\ff_1\ff_2^{-1})))(2u)g_\Phi(u)\\ 
          &=&ikk'\left(\frac{\wp'(\ff_2^2-\ff_1^2)}{\ff_1\ff_2\ff_3^2}\right)(2u)g_\Phi(u)\\ 
           &=& -i\left(\frac{\wp'}{\ff_1\ff_2\ff_3^2}\right)(2u)g_\Phi(u), 
\end{eqnarray*} 
since $e_2-e_1=-1/kk'$. Finally recall that $\wp'=-2\ff_1\ff_2\ff_3$.\ $\Box$ 
 
Using the eveness/oddness properties of $\ff_j(2u)$ at appropriate points, elementary considerations reveal that $G$ enjoys symmetries about $H_{\frac{1}{4}}$, $H_\frac{3}{4}$, $V_\frac{1}{4}$ and $V_\frac{3}{4}$: 
 
{\bf Proposition 8.2} (i)
$\di{G(\frac{\om_1}{4}-\bar{u})=G(\frac{\om_1}{4}+u)}$.    

(ii)
$\di{G(\frac{3\om_1}{4}-\bar{u})=G(\frac{3\om_1}{4}+u)}$.   
 
(iii) $\di{G(\frac{\om_2}{4}+\bar{u})=G(\frac{\om_2}{4}+u)}$. 
 
(iv) $\di{G(\frac{3\om_2}{4}+\bar{u})=G(\frac{3\om_2}{4}+u)}$. 
 
We now show that the integral  density of Gaussian curvature concentrates on $V_\frac{1}{4}\cup V_\frac{3}{4}$, as $k\ar1$, and on $H_{\frac{1}{4}}\cup H_\frac{3}{4}$, as $k\ar0$. The next result 
 follows immediately from 8.1 and
$\ff_3^2=(k'\ff_2+ik\ff_1)(k'\ff_2-ik\ff_1)$.

{\bf Proposition 8.3} 
$$ 
G(u)=\frac{4|g_\Phi'(u)|^2}{(1+|g_\Phi(u)|^2)^2}=\frac{8}{k^2|\wp(2u)-e_1|+k'^2|\wp(2u)-e_2|+|\wp(2u)-e_3|}. 
$$ 
It clarifies the exposition at this point to introduce the following reparameterizations. 
 Let $\rho_1:\C\ar\C$, be given by $\rho_1(z)=\om_1 z$ and $\rho_2:\C\ar\C$, be given by $\rho_2(z)=|\om_2| z$, moreover let $\tau_2=\om_2/\om_1$ and $\tau_1=\om_1/|\om_2|$: $\rho_1$ induces a biholomorphism $\C/(\Z 1+\Z\tau_2)\ar\C/(\Z\om_1+\Z\om_2)$, and $\rho_2$ a biholomorphism $\C/(\Z\tau_1+\Z i)\ar\C/(\Z\om_1+\Z\om_2)$. 
 
Let $V_a':=\{a+iy\ ; 0\leq y<\tau_2\}$, and $H_a':=\{x+ia\ ;\ 0\leq x <\tau_1\}$. Moreover, let $G_1, G_2$, be the integral densities induced by $\phi\circ\rho_1$, and $\phi\circ\rho_2$, respectively, i.e. 
$ 
G_j(z)=|\om_j|^2G(\om_jz)$, $j=1,2$.

{\bf Corollary 8.4}\ {\it For $u\in V_\frac{1}{4}\cup V_\frac{3}{4}$,  
$$ 
G(u)=\frac{4}{k'^2(\wp(2u)-e_2)}. 
$$ 
Hence,  $G_1\ar\infty$, on $V_\frac{1}{4}'\cup V_\frac{3}{4}'$,  as $k\ar1$.} 
 
{\it Proof.}\ For such $u$, $\ff_1(2u)\in i\R$, $\ff_2(2u)\in\R$ and $\ff_3(2u)\in\R$. Therefore \begin{eqnarray*} 
G(u)&=&\frac{8}{-k^2\ff_1(2u)^2+k'^2\ff_2(2u)^2+\ff_3(2u)^2}\\ 
    &=&\frac{8}{2k'^2\ff_2(2u)^2}.  
\end{eqnarray*} 
Now observe that for $z\in V_\frac{1}{4}'\cup V_\frac{3}{4}'$,  
$$ 
4k^2K(k)^2\leq\frac{\om_1^2}{k'^2(\wp(2\om_1z)-e_2)}\leq 4K(k)^2.\ \Box 
$$ 
Similiarly, we obtain: 
 
{{\bf Corollary 8.5}\ {\it For $u\in H_{\frac{1}{4}}\cup H_\frac{3}{4}$,  
$$ 
G(u)=-\frac{4}{k^2(\wp(2u)-e_1)}. 
$$ 
Hence, $G_2(u)\ar\infty$ on $H_{\frac{1}{4}}'\cup H_\frac{3}{4}'$,  as $k\ar0$.}

{\bf Corollary 8.6}\ {\it (i) At the quarter-period points we have: 
 
$$ 
G(\frac{\om_1}{4})=\frac{4k}{k'},\ \ \ \ 
G(\frac{\om_2}{4})=\frac{4k'}{k}\ \ {\rm and} \ \ 
G(\frac{\om_3}{4})=\frac{4}{kk'},\ \ {\rm etc}. 
$$ 
(ii) $\om_3/4$ is a maximum of \ $G$,  $\om_1/4$ and $\om_2/4$ are saddle-points.} 
 
{\it Remark.} Clearly, the behaviour of $G$ at the other quarter-period points may  be 
deduced from symmetry.
 
Comparing  
$$ 
G_1(z)\leq\frac{8\om_1^2}{|\wp(2\om_1z)-e_1|}, 
$$ 
it can be shown that for $z\not\in V_\frac{1}{4}'\cup V_\frac{3}{4}'$, $G_1(z)\ar0$, as $k\ar1$. 
Similiarly, for $z\not\in H_{\frac{1}{4}}'\cup H_\frac{3}{4}'$, 
$G_2(z)\ar0$, as $k\ar0$. This can be seen by inspecting power series. 
  
{\it Remark.} As $k\ar1$, narrower and narrower bands around $V_\frac{1}{4}'$
and $V_\frac{3}{4}' $  are stretched by the Gauss map to cover almost all of
the sphere.    

The underlying geometrical behaviour of $\gamma_\phi:\C/(\Z\om_1+\Z\om_2)\ar
S^2\subset\R^3$, elucidates the above results:    

{\bf Theorem 8.7}\  {\it (i)As $u$ passes  from $\om_1/4$ to $\om_1/4+\om_2$
along $V_\frac{1}{4}$, the Gauss map $\gamma_\phi(u)$ winds once around the
great circle  on $S^2$ that projects to the $\un{e}_1$-axis under
stereographic projection (from $-\un{e}_3$). The analogue holds on
$V_\frac{3}{4}$.    

(ii) As $u$ passes from $\om_2/4$ to $\om_2/4+\om_1$ along
$H_\frac{1}{4}$, the Gauss map $\gamma_\phi(u)$ winds once around the great
circle  on $S^2$ that projects to the $\un{e}_2$-axis under stereographic
projection  (from $-\un{e}_3$). } 

   {\it Proof.} This follows easily from the behaviour of
$g_\Phi$, which is real on $V_\frac{1}{4}\cup V_\frac{3}{4}$, while purely
imaginary on $H_\frac{1}{4}$.\ $\Box$  

We spell this out in more detail: Passing from $\om_1/4$ to $\om_3/4$,
$\phi$ maps onto the Higgs axis between $\beta_1$ and $\beta_3$, while the
normal vector turns through $90^o$. Between $\om_3/4$ and $\om_1/4+\om_2/2$,
$\phi$ maps back along the line segment, but as part of an intersecting
`sheet', see Figure 2. This is reflected by  the fact that the Gauss map turns
through another $90^o$, and arrives back at $\beta_2$ in the antipodal
direction. Now, this is all repeated between $\om_1/4+\om_2/2$ and
$\om_1/4+\om_2$, except that the values of the Gauss map are antipodal:
$\gamma_\phi(u+\om_2/2)=\alpha(\gamma_\phi(u))$, along $V_\frac{1}{4}$. This
reflects the fact that on the upper 1/2-rectangle of the period domain, $\phi$
maps onto the same surface in $\R^3$, but with opposite orientation. Of
course, the analogue holds along $V_\frac{3}{4}$. 

   As $k\ar1$, the length of
the line segment between $\beta_1$ and $\beta_3$ goes to zero, but the Gauss
map  makes one revolution along it for any $k\in(0,1)$. Thus as $k\ar1$, it
changes rapidly along the line segment, while as $k\ar0$, it winds slowly,
away from $\beta_3$.  

  Analogous remarks apply to the behaviour of the Gauss
map around $\Gamma_{\rm Star}(k)$. As $k\ar0$, $\Gamma_{\rm Star}(k)$ shrinks to
$0$, but the Gauss map winds around once on $\Gamma_{\rm Star}(k)$ for all
$k\in(0,1)$. See Figure 1. As $k\ar1$, $\Gamma_{\rm Star}(k)$ becomes very large.
(Note that $\gamma_\phi(u+\om_1/2)=\alpha(\gamma_\phi(u))$, on
$H_\frac{1}{4}$.)   

 Notice that the branch points $\pm\beta_3$ play a pivotal
role, connecting $\Gamma_{\rm Higgs}(k)$  to $\Gamma_{\rm Star}(k)$. Moreover the
integral density is large at these points in both limits.   

  We close this
section with some observations about the `bare' curvature ${\cal K}$ at
$\beta_1$ and $\beta_2$.  As $k\ar1$, $\beta_1$ accompanies $\beta_3$ to
infinity along the Higgs axis, while $\beta_2$, goes off to infinity along
$\un{e}_3$, and ${\cal K}(\om_2/4)$,  becomes attenuated. As $k\ar0$,
$\beta_1$ and $\beta_2$ `exchange roles'. More precisely, we have:   

 {\bf Proposition 8.8}\ {\it (i) There exist constants $\alpha_1, \alpha_2,
\alpha_3\in\R$, such that as $k\ar0$,   $$  |{\cal
K}(\frac{\om_1}{4}+\frac{h}{2})|\leq \alpha_1 k^5|h|^{-2}+\alpha_2
k^4+\alpha_3 k^3|h|^2+{\cal O}(|h|^4).  $$  (ii) There exist constants
$\beta_1, \beta_2, \beta_3\in\R$, such that as $k\ar1$,  $$  |{\cal
K}(\frac{\om_2}{4}+\frac{h}{2})|\leq \beta_1 k'^5|h|^{-2}+\beta_2 k'^4+\beta_3
k'^3|h|^2+{\cal O}(|h|^4).  $$}    It is not an optimal statement. The proof
follows from inspection of Laurent series.


\newpage


{\bf \S9. Scattering.}

The above results elucidate the behaviour of the  family of minimal surfaces
generated by the geodesic $C_1$ on ${\cal M}_2^0$, which corresponds to a
direct collision between two monopole particles, see \cite{ah} for further
details.  We will discuss the details elsewhere, however, we
make two comments: 

$\bullet$ Note that the $90^o$ scattering of the monopoles is reflected in the
behaviour of this family.

$\bullet$ The family includes a point (at the origin) at $k=0$. As $k\ar0$,
$\Gamma_{\rm Star}(k)$ shrinks to zero. This is better understood in terms of
the two curves on the null curve in $\C^3$ which project to $\Gamma_{\rm
Star}(k)$ .   These shrink to the two points $(\pm i\pi/4, 0, 0)$, (in
$0$-monopole coordinates) and the Gaussian curvature on the null curve
concentrates at these points.  Recall that the spectral curve of the axially
symmetric monopole gives the two stars of spectral affine null planes through
these points. Because  $(\pm i\pi/4, 0, 0)$ are not in $\R^3$, the
corresponding configurations of lines in $\R^3$ do not converge to the  star
through  the origin, see  \cite{sm5}.

{\it Remark.}  In fact it
is likely that many of the results described here should be understood
directly in terms of the null curve in $\C^3$. Cf. \cite{hurt1}, \cite{hurt3},
\cite{hurt4}, \cite{sm3.5} and \cite{sm4}.


\vspace{0.3in}

{\bf \S10. Charge 3.}

 Generically,
the spectral curve $S$, of a charge 3 monopole, is a smooth curve on ${\cal
C}(Q)\subset\p_3$ of genus 4. It is a canonical curve in $\p_3$, has degree 6
and  is a complete intersection of ${\cal C}(Q)$ with a cubic surface.

The  covering map
to $\p_1$ has six branch points. These come in antipodal pairs, interchanged
by $\tau$, and the three correponding spectral lines through    the origin in
$\R^3$ are perpendicular to the ends of the auxiliary minimal surface. (In
general, the auxiliary minimal surface of the generic charge $k$ monopole has
$k$ ends.)

Since $S$ is canonical, the points of hyperosculation are the
Weierstrass points on $S$.  The total weight of the Weierstrass points is 60.
Each of the six branch points of the covering map has weight 4, leaving 36 of
weight 1. $S$ is $\tau$-invariant and hence this gives 18  branch points on the
auxiliary minimal surface in $\R^3$. $18=6\times3$, where we group the
Weierstrass points in the same fibre. In summary we have:

{\bf Proposition 10.1}\  {\it The auxiliary minimal surface generated by a
generic \newline  monopole of charge 3 has three ends, total curvature
$-12\pi$, and 18 branch points in the metric.}

Of course, this does not
take us very far. Perhaps the next step is to consider the questions: 
  
$\bullet$ How does the area measure induced on the spectral curve by the Gauss
map relate to the configuration of Weierstrass points?

$\bullet$ How is
this reflected in the behaviour of the branch points and Gaussian curvature on
the surface in $\R^3$? 

However, perhaps all of this is better understood
in terms of the role played by the null curve in the geometry of the monopole's
complexification?



 \vspace{0.3in}  
{\bf Appendix} 
   \vspace{0.1in} 
  
\renewcommand{\arraystretch}{2.2}  \begin{tabular}{|l||c|c|c|}\hline 
\multicolumn{4}{|c|}{\bfseries Quarter-period Values: $g_2,\ g_3$ as in (8).}\\
\hline  \rule[-3mm]{0mm}{6mm}& $\di u=\frac{\om_1}{4}$ & $\di
u=\frac{\om_2}{4}$ & $\di u=\frac{\om_3}{4}$\\ \hline 
\rule[-3mm]{0mm}{9mm}$\di \wp(u)$ & $\di \frac{1+3k'+k'^2}{3kk'}$ & $\di
-\frac{1+3k+k^2}{3kk'}$ & $\di \frac{k^2-k'^2}{3kk'}-i$ \\ \hline 
\rule[-3mm]{0mm}{9mm}$\di \wp'(u)$ & $\di -\frac{2(1+k')}{k\sqrt{kk'}}$ & $\di
-\frac{2i(1+k)}{k'\sqrt{kk'}}$ & $\di \frac{2(k'+ik)}{\sqrt{kk'}}$\\ \hline 
$\di \wp''(u)$ & $\di \frac{4(1+k')}{k'(1-k')}$ & $\di \frac{4(1+k)}{k(1-k)}$
& $\di -8+\frac{4i(k'^2-k^2)}{kk'}$\\ \hline  \rule[-3mm]{0mm}{12mm}$\di
\wp'''(u)$ & $\di -\frac{8(1+3k'+k'^2)}{k'(1-k')\sqrt{kk'}}$ & $\di
\frac{8i(1+3k+k^2)}{k(1-k)\sqrt{kk'}}$ & $\di
8\left(\frac{5k^2-1}{k\sqrt{kk'}}-i\frac{5k'^2-1}{k'\sqrt{kk'}}\right)$\\
\hline  \end{tabular}

\vspace{0.2in}

{\it Remarks on the Figures.}

These were drawn with Mathematica. It is easy to write a short program using
(1)-(3), (7)-(9) and  (14). It is  instructive to do this and
experiment with the parameter values. We show the surfaces near parts of 
$\Gamma_{S_k}$ only. One can study larger regions but the pictures become
very complex, particularly when $k$ is not close to 0. It is also instructive
to plot $G$ and ${\cal K}$ for various values of $k$.

   \vspace{0.1in}     {\small    {\bf
Acknowledgements}    The author is very grateful to Patrick McCarthy for many
helpful conversations. He also thanks Werner Nahm and Paul Sutcliffe for
useful conversations. He is indebted to the Mathematics Faculty of the
University of Southampton for  their generous  hospitality and thanks Victor
Snaith for his invitation. He also thanks the ICTP, Trieste for their
hospitality during February and March of 2002.   }

\begin{figure}
\scalebox{0.6}{\includegraphics{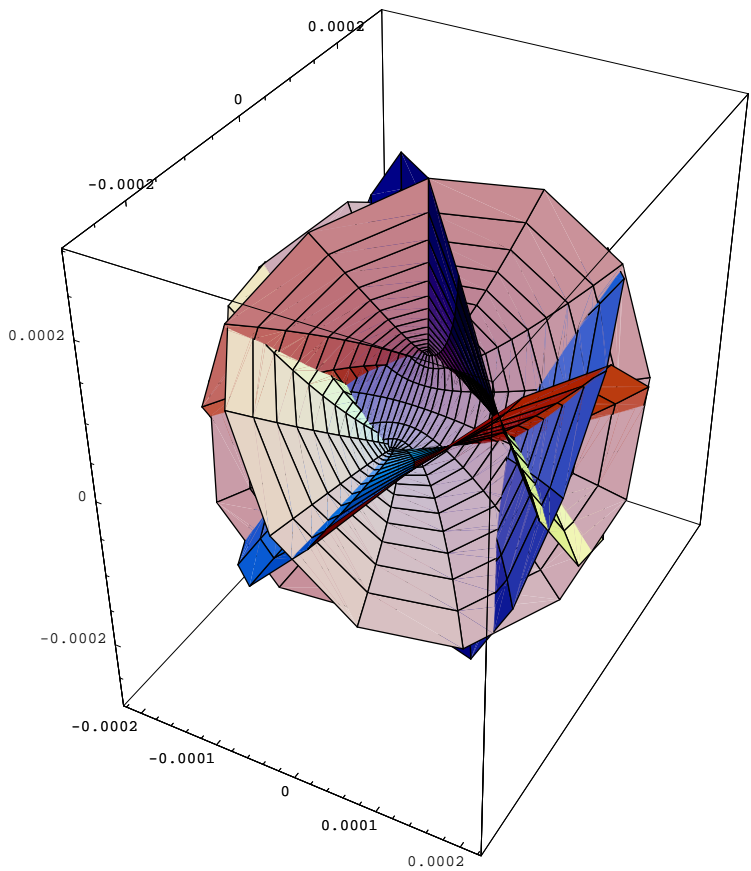}}
\caption{$k=0.01$,\ points close to $\Gamma_{\rm Star}(k)$: $0\leq x\leq\om_1,
\ \ \om_2/4-0.05\leq y\leq\om_2/4+0.05$} 
\end{figure}

\begin{figure}
\scalebox{0.6}{\includegraphics{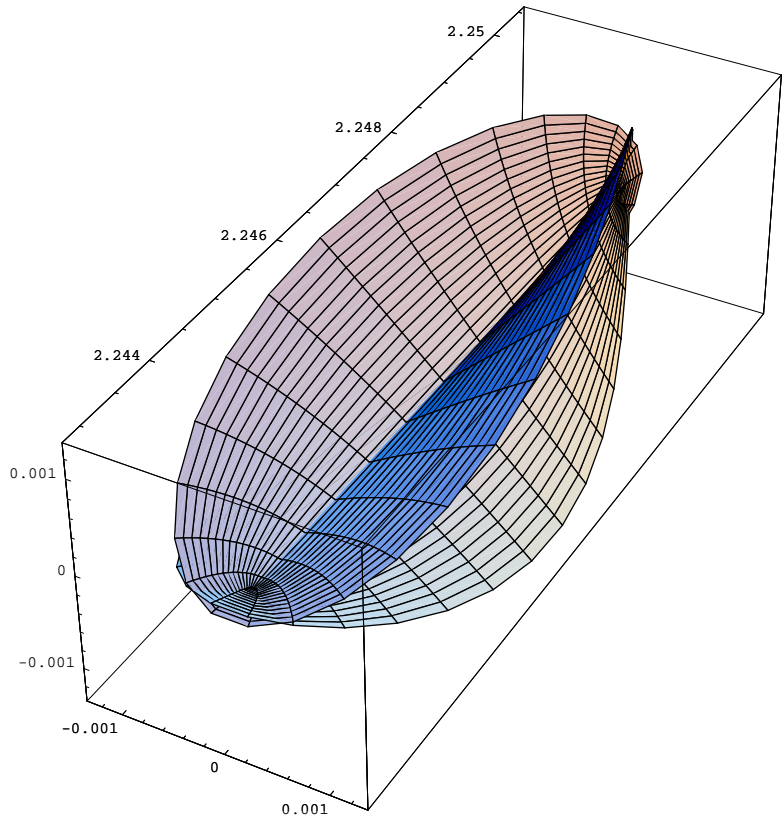}}
\caption{$k=0.999$,\ points close to (one half of) $\Gamma_{\rm Higgs}(k)$:
$\om_1/4-0.025\leq x\leq\om_1/4+0.025, \ \ 0\leq y\leq\om_2/2$} 
\end{figure}

\end{document}